\documentclass{elsarticle}
\usepackage{amsfonts}
\usepackage{graphicx}
\usepackage{amsmath}
\usepackage{amssymb}
\usepackage{amsthm}
\usepackage{enumerate}
\usepackage[all, cmtip]{xy}
\newtheorem{theorem}{Theorem}

\newtheorem{corollary}[theorem]{Corollary}

\newtheorem{example}[theorem]{Example}

\newtheorem{lemma}[theorem]{Lemma}

\newtheorem{proposition}[theorem]{Proposition}
\newtheorem{remark}[theorem]{Remark}

\newcommand{\End}{\operatorname{End}}
\begin{document}
\title{Rings and modules which are stable under automorphisms of their injective hulls}
\author[ne]{Noyan Er \corref{cor1}}
\ead{noyaner@yahoo.com}

\author[ss]{Surjeet Singh}
\ead{ossinghpal@yahoo.co.in}

\author[as]{Ashish K. Srivastava}
\ead{asrivas3@slu.edu}

 \cortext[cor1]{Corresponding author}

\address[ne]{Department of Mathematics,
University of Rio Grande, Rio Grande, OH 45674, USA.}
\address[ss]{House No. 424, Sector No. 35, A. Chandigarh 160036, India.}
\address[as]{Department of Mathematics and Computer Science, St Louis University, St Louis, MO 63103, USA.}

 \begin{abstract}
It is proved, among other results, that a prime right nonsingular
ring (in particular, a simple ring) $R$ is right self-injective if
$R_R$ is invariant under automorphisms of its injective hull. This
answers two questions raised by Singh and Srivastava, and Clark and
Huynh. An example is given to show that this  conclusion no longer
holds when prime ring is replaced by semiprime ring in the above
assumption. Also shown is that automorphism-invariant modules are
precisely pseudo-injective modules, answering a recent question of
Lee and Zhou. Furthermore, rings whose cyclic modules are
automorphism-invariant are investigated.
\end{abstract}
 \maketitle

\section{Introduction and Preliminaries}

\noindent Throughout,  $R$ will denote an associative ring with
identity and modules will be right modules. In \cite{df} Dickson and
Fuller studied modules which are invariant under automorphisms of
their injective hulls, when the underlying ring is a finite
dimensional algebra over a field with more than two elements. Such
modules over arbitrary rings were discussed by Lee and Zhou in
\cite{lz}, where they were called automorphism-invariant modules.
Thus, a module $M$ is called an automorphism-invariant module if $M$
is invariant under any automorphism of its injective hull. Clearly
every (quasi-)injective module is automorphism-invariant.

Dickson and Fuller had shown that if $R$ is a finite-dimensional
algebra over a field with more than two elements, then $R$ is of
right invariant module type if and only if every indecomposable
right $R$-module is automorphism-invariant. Recently, Singh and
Srivastava have investigated in \cite{AS} rings whose finitely
generated indecomposable right modules are automorphism-invariant,
and completely characterized indecomposable right Artinian rings
with this property. The dual notion of these modules has been
proposed by Singh and Srivastava in \cite{SS}.

The following questions are posed in the papers by Lee and Zhou
(\cite{lz}), Clark and Huynh (\cite{CH1}), and Singh and
Srivastava (\cite{AS}), respectively:\\

(Q1)  Is a simple ring $R$ such that $R_R$ is pseudo-injective right
self-injective \cite{CH1}?\\

(Q2)  Is a simple ring $R$ such that $R_R$ is automorphism-invariant
right self-injective \cite{AS}?\\

(Q3)  What is the structure of rings whose cyclic right modules are
automorphism-invariant \cite{AS}?\\

\noindent A module $M$ is called pseudo-injective if, for any
submodule $A$ of $M$, every monomorphism $A\rightarrow M$ can be
extended to some element of $\End(M)$. Pseudo-injective modules and
rings have been discussed by various authors (see, for example
\cite{AEJ}, \cite{hai}, \cite{JS}, \cite{MT}). Lee and Zhou showed
that a module $M$ is automorphism-invariant if and only if every
isomorphism between any two essential submodules of $M$ extends to
an automorphism of $M$ \cite{lz}. Thus it follows that
pseudo-injective modules are automorphism-invariant. Lee
and Zhou ask in \cite{lz} if the converse holds:\\

(Q4) Is an automorphism-invariant module pseudo-injective \cite{lz}?\\

\noindent In this paper, after proving a useful decomposition
theorem for an arbitrary automorphism-invariant module, we show that
a prime right nonsingular right automorphism-invariant ring is right
self-injective. Using this and the decomposition theorem, we
affirmatively answer the questions (Q1), (Q2) and (Q4). Also
obtained is a partial answer to (Q3).

For a property $P$ of modules, $R$ is said to have (or be) right $P$
if $R_R$ is a module with $P$. A closed submodule of a module $M$ is
one with no proper essential extensions in $M$. For submodules $A$
and $B$ of $M$, $B$ is said to be a complement of $A$ in $M$ if it
is maximal among submodules of $M$ trivially intersecting with $A$.
Complement submodules and closed submodules of $M$ coincide, and
being a closed submodule is a transitive property. An essential
closure of a submodule $A$ of a module $M$ is any closed submodule
of $M$ essentially containing $A$. In a nonsingular module, every
submodule has a unique essential closure. A module is called
square-free if it does not contain a direct sum of two nonzero
isomorphic submodules. Two modules are said to be orthogonal to each
other if they do not contain nonzero isomorphic submodules. A module
$B$ is said to be $A$-injective if every homomorphism from any
submodule $A'$ of $A$ into $B$ can be extended to an element of
$Hom(A,B)$. A detailed treatment of the above concepts and other
related facts can be found in \cite{ext} and \cite{mm}. Throughout
the paper, for a module $M$, $E(M)$ will denote the injective hull
of $M$.

\section{A decomposition theorem for automorphism-invariant modules}

Before proving our first main result, we will first give some useful
lemmas.

\begin{lemma}\label{l1} $($\cite[Lemma 7]{AS}$)$
If $M$ is an automorphism-invariant module with injective hull
$E(M)=E_1\oplus E_2\oplus E_3$ where $E_1\cong E_2$, then $M=(M\cap
E_1)\oplus (M\cap E_2)\oplus (M\cap E_3)$.
\end{lemma}

Lee and Zhou showed in \cite{lz} that whenever an
automorphism-invariant module $M$ has a decomposition $M=A\oplus B$,
$A$ and $B$ are relatively injective. This can be extended as
follows:

\begin{lemma}\label{zh2} If $M$ is an automorphism-invariant module and
$A$ and $B$ are closed submodules of $M$ with $A\cap B=0$, then $A$
and $B$ are relatively injective. Furthermore, for any monomorphism
$h:A\rightarrow M$ with $A\cap h(A)=0$, $h(A)$ is closed in $M$.
\end{lemma}

\emph{Proof.} First, let $K$ and $T$ be complements of each other in
$M$. Then, $E(M)=E_1\oplus E_2$, where $E_1=E(K)$ and $E_2=E(T)$.
Now let $f:E_1\rightarrow E_2$ be any homomorphism. Then the map
$g:E(M)\rightarrow E(M)$ defined by $g(x_1+x_2)=x_1+x_2+f(x_1)$
($x_i\in E_i$) is an automorphism, so that $f(K)=(g-1_{E(M)})(K)
\subseteq M$. Hence, $f(K)\subseteq E_2\cap M=T$. Therefore $T$ is
$K$-injective.

Now if $A$ and $B$ are closed submodules with zero intersection,
then, by the above argument, $A$ is injective relative to any
complement $C$ of $A$ containing $B$. Therefore, $A$ is
$B$-injective.

Finally, let $h:A\rightarrow M$ be a monomorphism with $h(A)\cap
A=0$, and pick any essential closure $K$ of $h(A)$. Since $A$ is
$K$-injective by the above arguments, $h^{-1}:h(A)\rightarrow A$
extends to a monomorphism $t:K\rightarrow A$. Therefore, we must
have $h(A)=K$.

\begin{theorem}\label{main} Let $M$ be an automorphism-invariant module.
Then the following hold:
\begin{itemize}
\item[(i)] $M=X \oplus Y$ where $X$ is quasi-injective and $Y$ is a
square-free module which is orthogonal to $X$. In this case, $X$ and
$Y$ are relatively injective modules.
\item[(ii)] If $M$ is nonsingular, then for any two
submodules $D_1$ and $D_2$ of $Y$ with $D_1\cap D_2=0$,
$Hom(D_1,D_2)=0$.
\item[(iii)] If $M$ is nonsingular, $Hom(X,Y)=0=Hom(Y,X)$.
\end{itemize}
\end{theorem}

\emph{Proof.} $(i)$ Let $\Gamma=\{(A,B,f): A,B\leq M, A\cap B=0,$
and $f:A\rightarrow B$ is an isomorphism$\}$.  Order $\Gamma$ as
follows: $(A,B,f)\leq (A',B',f')$ if $A\subseteq A'$, $B\subseteq
B'$, and $f'$ extends $f$. Then $\Gamma $ is inductive and there is
a maximal element in it, say $(A,B,f)$. Let $C'$ be a complement of
$A\oplus B$ in $M$. $C'$ must be square-free: Otherwise, there would
be nonzero submodules $X$ and $Y$ of $C'$ with $X\cap Y=0$, and an
isomorphism $\phi:X\rightarrow Y$. But then, $(A\oplus X, B\oplus Y,
f\oplus \phi)$ would contradict the maximality of $(A,B,f)$. So $C'$
is square-free. Now define $g:A\oplus B \oplus C' \rightarrow
A\oplus B\oplus C'$ via $g(a+b+c)=f^{-1}(b)+f(a)+c$ ($a\in A$, $b\in
B$, $c\in C'$). Since $M$ is automorphism-invariant, any isomorphism
between two essential submodules of $M$ extends to an automorphism
of $M$, whence $g$ extends to an automorphism $g'$ of $M$. Let $A'$
be a closed submodule of $M$ essentially containing $A$. If $A$ were
properly contained in $A'$, $g'_{|A'}$ would contradict the
maximality mentioned above. Thus, $A$ must be a closed submodule of
$M$. Since closed submodules are preserved under automorphisms, $B$
too is closed in $M$. Thus, by Lemma \ref{l1}, $M=(E(A)\cap M)
\oplus (E(B)\cap M) \oplus (E(C')\cap M)$. Then, $M=A \oplus B\oplus
C'$. Since direct summands of an automorphism-invariant module are
again automorphism-invariant, $A\oplus B$ is automorphism-invariant.
Now, by Lemma \ref{zh2}, it follows that $A$ and $B$ are relatively
injective. Since $A\cong B$, $A\oplus B$ is then quasi-injective.
Also, $A\oplus B$ and $C'$ are relatively injective modules. Next,
in a similar way to the above argument, one can find a maximal
monomorphism $t:B'\rightarrow B$ from a submodule $B' \subseteq C'$
into $B$. Since $B$ is $C'$-injective, $t$ can be monomorphically
extended to a closed submodule of $C'$ essentially containing $B'$.
By the maximality of $t$, this implies that $B'$ is closed in $C'$.
Also since $C'$ is $B$-injective, $t^{-1}$ extends monomorphically
to an essential closure, say $D$, of $t(B')$. Since $B'$ would then
be essential in the image of $D$, this implies that $t(B')$ is
closed in $B$. So $t(B')$ is a direct summand of $B$, since $B$ is
quasi-injective. And since $B$ is $C'$-injective, $t(B')$ is
$C'$-injective, hence $B'$ is a $C'$-injective submodule of $C'$.
Thus, $C'=B'\oplus C$ for some $C$. Now, we will show that $C$ and
$B$ are orthogonal: Assume that $C$ and $B$ have nonzero isomorphic
submodules $C_1$ and $B_1$. Then, by square-freeness of $C'$, $C_1$
and $B'$ are orthogonal modules, and thus, so are $B_1$ and $t(B')$,
so that we would have $B_1\cap t(B')=0$. This would contradict the
maximality of the monomorphism $t$. So $C$ and $B$ are orthogonal,
whence $C$ and $A\oplus B\oplus B'$ are orthogonal. Furthermore,
$A\oplus B\oplus B'$ is quasi-injective. Taking $X=A\oplus B\oplus
B'$ and $Y=C$, we obtain the desired conclusion.

$(ii)$ Let $f:D_1\rightarrow D_2$ be a nonzero homomorphism. By the
nonsingularity, $Ker(f)$ is closed in $D_1$ and there is some
submodule $L\neq 0$ of $D_1$ with $Ker(f)\cap L=0$. But then,
$L\cong f(L)\subseteq D_2$, contradicting the square-freeness of
$Y$. Now the conclusion follows.

$(iii)$ Similar to $(ii)$.

\begin{corollary} Any square-full automorphism-invariant module is
quasi-injective.
\end{corollary}

\begin{remark} \emph{Before the next result, note that in the proof of
Theorem \ref{main} $(ii)$, we have not used the assumption that $M$
is automorphism-invariant, so the statement holds for any
nonsingular square-free module.}
\end{remark}

Recall that a submodule $N$ of a module $M$ is called a fully
invariant submodule if, for every endomorphism $f$ of $M$,
$f(N)\subseteq N$.

\begin{theorem}\label{p} The following hold for a nonsingular square-free module $M$:
\begin{itemize}
\item[(i)] Every closed submodule of $M$ is a fully invariant submodule of $M$.
\item[(ii] If $M$ is automorphism-invariant, then for any family
$\{K_i:i\in I\}$ of closed submodules of $M$ (not necessarily
independent), the submodule $\Sigma_{i\in I}K_i$ is
automorphism-invariant.
\end{itemize}
\end{theorem}

\emph{Proof.} First, assume that $M$ is square-free and nonsingular.
Let $K$ be a closed submodule of $M$ and $T$ be a complement in $M$
of $K$. Suppose that $f\in \End(M)$ with $f(K)\nsubseteq K$. Let
$\pi :E(M) \rightarrow E(T)$ be the obvious projection with
$Ker(\pi)=E(K)$. Since $K$ is not essential in $f(K)+K$, we have
$\pi(K+f(K))\neq 0$, implying that $\pi(f(K))\neq 0$, whence $N=T
\cap \pi(f(K))\neq 0$. Then, for $N'=\{x\in K: \pi f(x)\in T\}$, we
have $Hom(N',N)\neq 0$, contradicting the assertion preceding this
theorem. This proves $(i)$.

Now assume, furthermore, that $M$ is automorphism-invariant, and let
$\{K_i: i\in I\}$ be any family of closed submodules of $M$, and $g$
be an automorphism of $E(\Sigma_{i \in I}K_i)$. Clearly, $g$ can be
extended to an automorphism $g'$ of $E(M)$. Since $M$ is
automorphism-invariant, we have $g'(M)\subseteq M$. Then, by $(i)$,
$g(K_i)=g'(K_i)\subseteq K_i$ for all $i\in I$. This proves $(ii)$.

\section{Nonsingular automorphism-invariant rings}

In this section we will prove a theorem describing right nonsingular
automorphism-invariant rings and answer two questions raised by
Singh and Srivastava in \cite{AS}, and by Clark and Huynh in
\cite{CH1} concerning when an automorphism-invariant or a
pseudo-injective ring is self-injective.

\begin{theorem}\label{f}
If $R$ is a right nonsingular right automorphism-invariant ring,
then $R\cong S \times T$, where $S$ and $T$ are rings with the
following properties:
\begin{itemize}
\item[(i)] $S$ is a right self-injective ring,
\item[(ii)] $T_T$ is square-free, and
\item[(iii)] Any sum of closed right ideals of $T$ is a two
sided-ideal which is automorphism-invariant as a right $T$-module.
\item[(iv)] For any prime ideal $P$ of $T$ which is not essential in $T_T$,
$\frac{T}{P}$ is a division ring.
 \end{itemize}
\end{theorem}

\emph{Proof.}  By Theorem \ref{main}, $R= eR\oplus (1-e)R$ for some
idempotent $e\in R$, where
 $eR$ is quasi-injective, $(1-e)R$ is
square-free and \[Hom(eR,(1-e)R)=0=Hom((1-e)R,eR).\] Hence, $S=eR$
and $T=(1-e)R$ are ideals. Now we have $(i)$ and $(ii)$. Also,
$(iii)$ follows from Theorem \ref{p}.

We now prove $(iv)$: Let $P$ be a prime ideal of $T$ which is not
essential as a right ideal. Take a complement $N$ of $P$ in $T_T$.
If $N$ were not uniform, there would be two nonzero closed right
ideals in $N$, say $X$ and $Y$ with $X\cap Y=0$. They would then be
ideals by the above argument. But this would contradict the
primeness of $P$. So $N$ is a uniform right ideal of $T$. Also note
that $P$ is a closed submodule of $T_T$, because if $P'$ is any
essential extension of $P$, we have $P'N=0$, implying that $P'=P$.
So $P$ is closed in $T_T$, and hence it is a complement in $T_T$ of
$N$. Since $\frac{N\oplus P}{P}$ is essential in $\frac{T}{P}$, this
implies that the ring $\frac{T}{P}$ is right uniform. Furthermore,
$N$ is a nonsingular uniform automorphism-invariant $T$-module, so
that every nonzero homomorphism between any two submodules is an
isomorphism between essential submodules, and thus it extends to an
automorphism of $N$. Therefore, $N$ is a quasi-injective uniform
nonsingular $\frac{T}{P}$-module, and thus its endomorphism ring is
a division ring. Since $\frac{T}{P}$ essentially contains the
nonsingular right ideal $\frac{N\oplus P}{P}$, it is now a prime
right uniform and right nonsingular ring (hence a prime right Goldie
ring) with the quasi-injective essential right ideal $\frac{N\oplus
P}{P}$. But then $\frac{N\oplus P}{P}$ is injective, implying that
$P\oplus N=T$. In fact, since $\End_{\frac{T}{P}}(N)$ is a division
ring, $\frac{T}{P}$ is a division ring.  In particular, $N$ is a
simple right ideal and $P$ is a maximal right ideal of $T$.

\begin{theorem}\label{d} If $R$ is a prime right non-singular, right automorphism-invariant
ring, then $R$ is right self-injective.
\end{theorem}

\emph{Proof.} By Theorem \ref{f} and primeness, it suffices to look
at the case when $R_R$ is square-free: If $R_R$ were not uniform,
there would be two closed nonzero right ideals $A$ and $B$ with
$A\cap B=0$. But then $A$ and $B$ would be ideals, whence $AB=0$,
contradicting primeness. So $R_R$ is uniform, nonsingular and
automorphism invariant. Now it follows, in the same way as in the
proof of Theorem \ref{f},
that $R$ is right self-injective.\\

 The following
example shows that the conclusion of Theorem \ref{d} fails if we
take a semiprime ring instead of a prime one.

\begin{example} \emph{Let $S=\prod_{n\in \mathbb{N}}\mathbb{Z}_2$, and
$R=\{(x_n)_{n \in \mathbb{N}}:$ all except finitely many $x_n$ are
equal to some $a\in \mathbb{Z}_2 $ $\}$. Then $S$ is a commutative
self-injective ring with $S=E(R_R)$ with only one automorphism,
namely the identity. Thus, $R$ is an automorphism-invariant,
semiprime nonsingular ring, but it is not self-injective. Teply
constructed in \cite{MT} the first example of a pseudo-injective
module which is not quasi-injective. In fact, the ring $R$ here is a
new example of pseudo-injective ring which is not self-injective, by
Theorem \ref{pseudo} below.}
\end{example}

The following corollary answers the question of Singh and Srivastava
in \cite{AS}.

\begin{corollary} A simple right automorphism-invariant ring is
right self-injective.
\end{corollary}

 The next corollary answers the question raised by Clark and
Huynh in \cite[Remark 3.4]{CH1}.

\begin{corollary} A simple right pseudo-injective ring is right self-injective.\end{corollary}

\section{Rings whose cyclic modules are automorphism-invariant}

Characterizing rings via homological properties of their cyclic
modules is a problem that has been studied extensively in the last
fifty years. A most recent account of results related to this
prototypical problem may be found in \cite{JST}, and a recent
addition in \cite{glas}. Another question raised in \cite{AS} is the
following: What is the structure of rings whose cyclic right modules
are automorphism-invariant? The next result addresses this question.

\begin{theorem}\label{cyc} Let $R$ be a ring over which every cyclic right
$R$-module is automorphism-invariant. Then $R\cong S \times T$,
where $S$ is a semisimple artinian ring, and $T$ is a right
square-free ring such that, for any two closed right ideals $X$ and
$Y$ of $T$ with $X\cap Y=0$, $Hom(X,Y)=0$. In particular, all
idempotents of $T$ are central.
\end{theorem}

\emph{Proof.} By the proof of Theorem \ref{main}, we have a
decomposition $R_R= A\oplus B\oplus B'\oplus C$, where $A\cong B$,
$B'$ is isomorphic to a submodule of $B$,  and $C$ is square-free
and $A\oplus B\oplus B'$ and $C$ are orthogonal. Let $Z$ be a right
ideal in $A$. Then $\frac{R}{Z}\cong \frac{A}{Z} \oplus B \oplus
B'\oplus C$ is automorphism-invariant by assumption. Then, by Lemma
\ref{zh2}, $\frac{A}{Z}$ is $B$-injective, whence $A$-injective.
Similarly, all factors of $B$, $B'$, and $C$ are $A$-injective as
well.

Now, $A$ is a cyclic projective module all of whose factors are
$A$-injective (and in particular, quasi-injective). So, by
\cite[Corollary 9.3 (ii)]{ext}, $A=U_1\oplus ... \oplus U_n$, where
$U_i$ are uniform modules. Take an arbitrary nonzero cyclic
submodule $U$ of $U_i$, for any $i$. Since $U$ is a sum of factors
of $A$, $B$, $B'$ and $C$, it contains a nonzero factor of one of
them, call $U'$. By the above paragraph, $U'$ is $A$-injective, so
it splits in $U_i$. Thus, $U'=U=U_i$, showing that $U_i$ is simple,
whence $A\oplus B\oplus B'$ is semisimple. Since $A\oplus B\oplus
B'$ and $C$ are orthogonal projective modules and the former is now
semisimple, there are no nonzero homomorphisms between them.
Therefore, $A\oplus B\oplus B'$ and $C$ are ideals. So now we have
the ring direct sum $R=S\oplus T$ where $S=A\oplus B\oplus B'$ and
$T=C$.

Now let $X$ and $Y$ be closed right ideals of $T$ such that $X\cap
Y=0$, and let $f:X\rightarrow Y$ be any homomorphism. Set $Y'=f(X)$.
This induces an isomorphism $\overline{f}:\frac{X}{K}\rightarrow
Y'$, where $K=Ker(f)$. It is clear that $\frac{X}{K}$ is a closed
submodule of $\frac{T}{K}$. Also, since $T_T$ is square-free, $K$ is
essential in $X$. Choose a complement $\frac{U}{K}$ of
$\frac{X}{K}\oplus \frac{Y'\oplus K}{K}$ in $\frac{T}{K}$. Since
$\frac{T}{K}$ is automorphism-invariant by assumption and
$\frac{X}{K}\cong Y'\cong \frac{Y'\oplus K}{K}$, by the last part of
Lemma \ref{zh2}, $\frac{Y'\oplus K}{K}$ is closed in $\frac{T}{K}$.
Applying Lemma \ref{l1}, we obtain $\frac{T}{K}=\frac{X}{K}\oplus
\frac{Y'\oplus K}{K} \oplus \frac{U}{K}$. Since $Y'\cap
(X+U)\subseteq Y'\cap K=0$, we have $T= Y' \oplus (X+U)$. So $Y'_T$
is projective, whence the map $f$ above splits. However, since $K$
is essential in $X$, we have $f=0$. So, $Hom(X,Y)=0$. In particular,
if $T_T=X\oplus Y$, we have $XY=YX=0$,
whence $X$ and $Y$ are ideals.\\

Using an alternative argument to the one in the second paragraph of
the above proof, we can generalize the decomposition in the theorem
as follows:

\begin{proposition} Let $M$ be a module satisfying any one of the following conditions:
\begin{itemize}
\item[(i)] $M$ is cyclic with all factors automorphism-invariant, and generates its cyclic subfactors, or
\item[(ii)] $M$ is any automorphism-invariant module whose 2-generated subfactors are
automorphism-invariant.
\end{itemize}
Then $M=X\oplus Y$, where $X$ is semisimple, $Y$ is square-free, and
$X$ and $Y$ are orthogonal.
\end{proposition}

\emph{Proof.}  First note that, by the proof of Theorem \ref{main},
we have a decomposition $M=A\oplus B \oplus B'\oplus C$, where
$A\cong B$, $B'$ embeds in $B$, and $C$ is square-free and
orthogonal to $A\oplus B\oplus B'$.

$(i)$ In this case, in the same way as in the first paragraph of the
proof of Theorem \ref{cyc}, all factors of the modules $B$ ($\cong
A$), $B'$ and $C$ are $A$-injective. Now let $A'$ be any factor of
$A$ and $D$ be a cyclic submodule of $A'$. Since $D$ is generated by
$M$, $D=D_1 + ...+ D_n$, where each $D_i$ is a factor of $B$, $B'$
or $C$. Since $D_1$ is $A$-injective (whence $A'$-injective),
$D_1\oplus D_1'=A'$ for some submodule $D_1'$ of $A'$. Letting
$\pi:D_1\oplus D_1' \rightarrow D_1'$ be the obvious projection, we
have $D=D_1 \oplus (\pi(D_2) + ... + \pi(D_n))$. Each $\pi(D_k)$
again being a factor of $B$, $B'$ of $C$, it is $A$-injective,
whence $D_1'$-injective. By induction on $n$, we obtain that $D$ is
a direct sum of $A$-injective cyclic modules. Then $D$ is
$A$-injective. Now we have shown that each cyclic subfactor of $A$
is $A$-injective. By \cite[Corollary 7.14]{ext}, $A$ is semisimple.
Therefore, $A\oplus B\oplus B'$ is semisimple, as well. Now set
$X=A\oplus B\oplus B'$ and $Y=C$.

$(ii)$ Let $D\subseteq L$ be submodules of $A$ with $\frac{L}{D}$
cyclic, and $T$ be a cyclic submodule of $B$. By assumption,
$\frac{L}{D}\oplus T$ is automorphism-invariant, whence
$\frac{L}{D}$ is $T$-injective. Then, cyclic subfactors of $A$ are
 $B$-injective, hence $A$-injective. Again, by \cite[Corollary
7.14]{ext}, $A$ is semisimple. The conclusion follows in the same
way as above.\\

\section{Pseudo-injective modules and automorphism invariant modules coincide}

In \cite{lz} Lee and Zhou raise the following question: Is an
automorphism-invariant module pseudo-injective? In the next theorem,
this question is answered affirmatively, also settling
\cite[Question 2]{AS}.

First we recall a useful lemma.

\begin{lemma}\label{relinj}(\cite[Lemma 7.5]{ext}) Let $M=A\oplus B$. Then $A$ is $B$-injective if and only if
for any submodule $C$ of $M$ with $A\cap C=0$, there exists some
submodule $D$ of $M$ such that $C\subseteq D$ and $A\oplus D=M$.
\end{lemma}

\begin{lemma}\label{A} Assume that $M=A\oplus B$ where, $A$ and $B$ are orthogonal to each other.
For any submodule $C$ of $M$ and any monomorphism $f:C\rightarrow M$
the following assertions hold:
\begin{itemize}
\item[(i)] $f(C\cap B)\cap B$ is essential in $f(C\cap B)$.
\item[(ii)] If $B$ is square-free, then $f(C\cap B)\cap (C\cap B)$ is
essential in both $f(C\cap B)$ and $C\cap B$.
\end{itemize}
\end{lemma}

\emph{Proof.} Let $C$ be a submodule of $M$ and $f:C\rightarrow M$
be a monomorphism. Assume $D$ is a submodule of $f(C\cap B)$ with
$D\cap B=0$. Then $D$ is embedded (via the obvious projection
$A\oplus B \rightarrow A$) into $A$. But $D$ is also isomorphic to a
submodule of $C\cap B$. This implies, by orthogonality, that $D=0$.
This proves $(i)$.

Now, assume $X$ is a nonzero submodule of $f(C\cap B)$ with $X\cap
(C\cap B)=0$. Then by $(i)$, $X\cap B\neq 0$, and now $(X\cap B)^2$
embeds in $(X\cap B)\oplus(C\cap B)\subseteq B$, a contradiction to
the assumption that $B$ is square-free. Hence, $f(C\cap B)\cap
(C\cap B)$ is essential in $f(C\cap B)$. One can see similarly that
$f(C\cap B)\cap (C\cap B)$ is also essential in $C\cap B$. This
proves $(ii)$.

\begin{theorem}\label{pseudo}
A module $M$ is automorphism-invariant if and only if it is
pseudo-injective.
\end{theorem}

\emph{Proof.} The fact that pseudo-injective modules are
automorphism-invariant follows from \cite{lz}. So, let $M$ be
automorphism-invariant, $C$ be a submodule of $M$, and
$f:C\rightarrow M$ be a monomorphism. By Theorem \ref{main},
$M=A\oplus B$, where $A$ is quasi-injective, $B$ is square-free
automorphism-invariant, and $A$ and $B$ are relatively injective.
Now let $K$ be a complement in $B$ of $f(C\cap B)\cap(C\cap B)$.
Then, by Lemma \ref{A} $(ii)$, $K\oplus [f(C\cap B)\cap(C\cap B)]$
is essential in both $K\oplus (C\cap B)$ (hence in $B$) and $K\oplus
f(C\cap B)$. This implies $[K\oplus f(C\cap B)]\cap A=0$, and
 $K\oplus f(C\cap B)\oplus A$ is essential in $M$.

Since $A$ is $B$-injective, then by Lemma \ref{relinj}, there exists
a submodule $B'$ of $M$ such that $f(C\cap B)\oplus K\subseteq B'$
and $M=A\oplus B'$. In this case, $B'\cong B$. By the above
paragraph, $f(C\cap B) \oplus K$ is essential in $B'$.
 Since $B$ is automorphism-invariant, the isomorphism $f_{|C\cap B}\oplus 1_K : (C\cap B)\oplus K\rightarrow
 f(C\cap B)\oplus K$ extends to some isomorphism $f':B\rightarrow B'$. So
 now $f'_{|C\cap B}=f_{|C\cap B}$.

The map $g:C+B\rightarrow f(C)+B'$ defined by $g(c+b)=f(c)+f'(b)$
($c\in C$, $b\in B$) is well-defined and extends $f$. Now let
$\pi:A\oplus B\rightarrow A$ be the obvious projection. Then
$B+C=B\oplus \pi(C)$. Note that $\pi(C)=(B+C)\cap A$. Since $A$ and
$B$ are both $A$-injective, then $M$ is $A$-injective, whence
$g_{|\pi(C)}:\pi(C)\rightarrow M$ extends to some $g':A\rightarrow
M$. Clearly, $g'_{|\pi(C)}=g'_{|(B+C)\cap A}=g_{|(B+C)\cap A}$. Now
we define $\psi:M\rightarrow M$ as follows: For $a\in A, x\in B+C$,
$\psi(a+x)=g'(a)+g(x)$. $\psi$ is the desired extension of $f$ to
$M$. Therefore $M$ is pseudo-injective.\\

Since pseudo-injective modules are known to satisfy the property
(C$_2$) by \cite{hai}, this also yields the affirmative answer to
another question in \cite{AS}.\\

\textbf{Acknowledgement} The authors would like to thank the
referee, whose careful reading and thoughtful comments have
helped improve the paper.\\

\end{document}